\documentclass[12pt]{amsart}

\usepackage{amsmath,amssymb,amsthm}
\usepackage[dvipdfm]{graphicx}
\usepackage{enumerate}
\usepackage[dvips]{color}
\usepackage{version}

\newtheorem{Thm}{Theorem}[section]
\newtheorem{Lem}[Thm]{Lemma}

\newtheorem{Prop}[Thm]{Proposition}

\newtheorem{Conj}[Thm]{Conjecture}
\newtheorem{Cond}[Thm]{Condition}

\theoremstyle{remark}
\newtheorem{Rem}[Thm]{Remark}
\newtheorem{Ex}[Thm]{Example}
\theoremstyle{definition}
\newtheorem{Def}[Thm]{Definition}

\newtheorem*{ack}{Acknowledgments}

\newcommand{\lct}{\mathop{\mathrm{lct}}\nolimits}
\newcommand{\ct}{\mathop{\mathrm{ct}}\nolimits}

\newcommand{\Exc}{\mathop{\mathrm{Exc}}\nolimits}
\newcommand{\Bir}{\mathop{\mathrm{Bir}}\nolimits}
\newcommand{\Aut}{\mathop{\mathrm{Aut}}\nolimits}

\newcommand{\Pic}{\mathop{\mathrm{Pic}}\nolimits}
\newcommand{\codim}{\mathop{\mathrm{codim}}\nolimits}

\newcommand{\Supp}{\mathop{\mathrm{Supp}}\nolimits}

\newcommand{\val}{\mathop{\mathrm{val}}\nolimits}

\newcommand{\DF}{\mathop{\mathrm{DF}}\nolimits}
\newcommand{\Sesh}{\mathop{\mathrm{Sesh}}\nolimits}

\excludeversion{Eliminates}

\excludeversion{N.B.}

\excludeversion{memo}

\begin{document}

\title[Birational rigidity and GIT]
{Birational superrigidity and \\ slope stability of Fano manifolds}

\author{Yuji Odaka}
\address{Department of Mathematics, Faculty of Science, Kyoto University}
\email{yodaka@math.kyoto-u.ac.jp}

\author{Takuzo Okada}
\address{Department of Mathematics, Faculty of Science and Engineering, Saga University} 
\email{okada@cc.saga-u.ac.jp}
\date{8th, September, 2011. Revised December, 2012.}
\maketitle

%%%abstract%%%
\begin{abstract}
We prove a relation between birational superrigidity of Fano manifold and its 
slope stability in the sense of Ross-Thomas \cite{RT07}. 
\end{abstract}
%%%%%%%%%%
\tableofcontents

\section{Introduction }
The concept of \textit{birational} (\textit{super})\textit{rigidity} of Fano manifolds (or of  
Mori fiber spaces, in general) was introduced to extend the work of Iskovskikh-Manin \cite{IM71} for quartic threefolds. The concept emerged in the study of 
the rationality problem for Fano manifolds by analyzing birational maps between such spaces. 

The purpose of this paper is to show a relation between birational 
(super)rigidity 
and GIT stability, which is unexpected because of the 
different nature of their origins. More precisely, in this paper we study 
\textit{slope stability} of polarized varieties, which was introduced by Ross-Thomas  (cf.\  \cite{RT07}) as an analogue of Mumford-Takemoto's slope stability of vector bundles. It is also a weaker version of K-stability, which was first formulated by Tian in \cite{Tia97} and later reformulated and generalized by Donaldson \cite{Don02}. 
(Ross-Thomas \cite{RT07} follow the formulation of \cite{Don02}.) 

Our main result is the following. 

\begin{Thm}\label{thm:main}
Let $X$ be a birationally superrigid Fano manifold of Fano index $1$. 
If $|-K_X|$ is base point free, then $(X,\mathcal{O}_{X}(-K_X))$ is slope stable. 
\end{Thm}
\noindent
We remark that the assumptions in Theorem \ref{thm:main} on the index and the base point freeness of the anticanonical linear system $|-K_X|$ seem to be mild.
As far as the authors know, every Fano manifold which has been proved to be birationally superrigid satisfies both assumptions (see section \ref{subsec:bir.rigid} for examples of birationally superrigid Fano manifolds). 

\color{black}{Actually we prove the following stronger but technical result from which Theorem \ref{thm:main} follows. 

\begin{Thm}\label{str.main}
Let $X$ be a Fano manifold of Picard rank $1$ and index $1$ with no log maximal singularity (see section \ref{subsec:bir.rigid} for the definition).
If $|-K_X|$ is base point free, then $(X, \mathcal{O}_X (-K_X))$ is slope stable.
\end{Thm}}\color{black}{}

Recall that the motivation for introducing K-stability is to formulate the 
following relation with the existence of K\"ahler metrics, which has been called 
the Yau-Tian-Donaldson conjecture. 
\begin{Conj}\label{fact:YTD}
An arbitrary Fano manifold $X$ with discrete automorphism group admits 
a K\"ahler-Einstein metric if and only if $(X,-K_X)$ 
is K-stable. 
\end{Conj}
\noindent
The only if direction was finally proved by \cite{Stp08} and 
proofs of the if direction were recently announced by \cite{CDS12}, \cite{Tia12}. 
Note that every birationally superrigid Fano manifold has discrete automorphism group  since it is not (birationally) ruled. 
We also refer the readers to \cite{Li13} for a differential geometric counterpart of $K$-semistability as well.

We remark that the following example shows that our result \textit{cannot} be derived from the main result of \cite{OS10}. 
\begin{Ex}
Let $X$ be a smooth projective hypersurface of dimension $n$ and degree $n+1$ in $\mathbb{P}^{n+1}$. 
After the papers \cite{IM71} and \cite{P98}, it has been completely proved in \cite{dF11} that $X$ is birationally superrigid for $n \ge 3$.
%Due to \cite{P98}, a general $X$ is birationally superrigid for $n \ge 4$ and this is conjectured to hold for every nonsingular $X$.
%In fact, every smooth $X$ is proved to be birationally superrigid for $n = 3$ by \cite{IM71} and for $4 \le n \le 12$ by \cite{dFEM03}. 
On the other hand, 
it is known that the alpha invariant $\alpha(X)$ of Tian (\cite{Tia87}, cf. also \cite{OS10})  is $\frac{n}{n+1}$ if $X$ contains some \textit{generalized Eckardt points} (or equivalently, hyperplane sections of cone type) so that strict stability does \textit{not} directly follow from 
\cite{OS10}. We are grateful to Professor Constantin Shramov for pointing this out to us. We refer the readers to \cite{OS10} for basics of the alpha invariant and its algebro-geometric  version. \end{Ex}

Our proof of Theorem \ref{thm:main} is similar to that of \cite{OS10}. 
Recall that the two fundamental observations in \cite{OS10} are that: 
\begin{itemize}
\item Certain explicit upper bounds for the Seshadri constants imply K-stability of a $\mathbb{Q}$-Fano variety (see \cite[Corollary 4.4]{OS10}). 

\item Mildness of singularities of pluri-anticanonical divisors gives upper bounds for 
the Seshadri constants. 
\end{itemize}
\noindent
We combine these observations to prove Theorem \ref{thm:main}, which is possible since birational superrigidity asserts 
certain mildness of singularities of pluri-anticanonical divisors as we will review in 
subsection \ref{b.rig.equivalence}. 

In the next section, we give some basic definitions and review the background. 
In section \ref{div.sec}, we prove stability along divisors and in section \ref{h.cod.sec}, we prove stability along higher codimensional loci. The last section proposes a more general conjecture about stability of Fano manifolds. 

We work over the field of complex numbers $\mathbb{C}$ throughout this paper. 
 
\begin{ack}
We are grateful to Professors Shigefumi Mori, 
Constantin Shramov and Alexander Pukhlikov for helpful comments. We would like to thank Doctors 
Jesus Martinez-Garcia and Julius Ross for their careful reading of the draft versions and 
Professor Ivan Cheltsov for his interest in this work. 
We also would like to thank the referee for her/his efforts to point out many errors.
The first author is partially supported by the Grant-in-Aid for Scientific Research (KAKENHI No.\ 21-3748)  and the Grant-in-Aid for JSPS fellows 
(DC1). The second author is partially supported by the Grand-in-Aid for Scientific Research (KAKENHI No.\  23-2053) and the Grand-in-Aid for JSPS fellows (PD). 

\end{ack}

\section{Preliminaries}\label{sec:Preliminary}

\subsection{Birational (super)rigidity}\label{subsec:bir.rigid}
In this subsection, we recall the definition of birational (super)rigidity and its basic properties.
We refer the readers to \cite{P07} for an introduction to this subject.

\begin{Def}
Let $\pi \colon V \to S$ be a projective surjective morphism between normal projective varieties with connected fibers.
We say that $\pi \colon V \to S$ is a {\it Mori fiber space} if
\begin{enumerate}
\item $V$ is $\mathbb{Q}$-factorial and has at most terminal singularities,
\item $-K_V$ is $\pi$-ample,
\item $\dim S < \dim V$, and
\item the relative Picard number $\rho (V/S)$ is $1$.
\end{enumerate}
\end{Def}

Throughout this subsection, let $X$ be a $\mathbb{Q}$-factorial $\mathbb{Q}$-Fano variety with Picard number $1$ and with at most terminal singularities.
Note that $X$, together with the structure morphism (to a point), 
can be seen as a Mori fiber space. 
Although birational (super)rigidity can be defined for any Mori fiber space, we only give the definition for $\mathbb{Q}$-Fano varieties with Picard number one. 

\begin{Def}
We say that $X$ is {\it birationally rigid} if for any birational map 
$\varphi \colon X \dashrightarrow X'$ to a Mori fiber space there is a birational self-map $\tau \colon X \dashrightarrow X$ such that $\varphi \circ \tau \colon X \dashrightarrow X'$ can be extended to an isomorphism.
We say that $X$ is {\it birationally superrigid} if in the above definition of birational rigidity one can always take $\tau = \operatorname{id}_X$.
\end{Def}

It is easy to see that $X$ is birationally superrigid if and only if $X$ is birationally rigid and $\Bir (X) = \Aut (X)$.
Let $\mathcal{H}$ be a movable linear system on $X$, that is, a linear system without fixed components. 
We define $\mu = \mu (X, \mathcal{H})$ to be the rational number for which $\mu K_X + \mathcal{H} \equiv 0$, where $\equiv$ denotes the numerical equivalence.
Let $\lambda$ be a nonnegative rational number.
We say that a pair $(X, \lambda \mathcal{H})$ is {\it terminal} (resp.\ {\it canonical}, resp.\ {\it log canonical}) if every rational number $a (X, \lambda \mathcal{H}, E)$ determined by 
\[
K_V + \lambda f_*^{-1} \mathcal{H} = f^*(K_X + \lambda \mathcal{H}) + \sum a (X, \lambda \mathcal{H}, E) E,
\]
is positive (resp.\ non-negative, \color{black}{resp.\ $\ge -1$}\color{black}{}) for every birational morphism $f \colon V \to X$, where $E$ runs over the $f$-exceptional prime divisors. 

The \textit{canonical threshold} (resp.\ \textit{log canonical threshold}) of the pair $(X,\mathcal{H})$ is defined to be the number 
$$
{\rm ct}(X,\mathcal{H}):=\sup\{\lambda\in \mathbb{Q}_{>0}\mid (X,\lambda \mathcal{H}) \mbox{ is canonical}\} 
$$
$$
(\text{resp.\ }{\rm lct}(X,\mathcal{H}):=\sup\{\lambda\in \mathbb{Q}_{>0}\mid (X,\lambda \mathcal{H}) \mbox{ is log canonical}\}). 
$$

\begin{Def}
We say that $X$ has a {\it maximal singularity} (resp.\ {\it log maximal singularity}) if there is a movable linear system $\mathcal{H}$ on $X$ such that $(X, \frac{1}{\mu} \mathcal{H})$ is not canonical \color{black}{(resp.\ not log canonical)}\color{black}{} for $\mu = \mu (X,\mathcal{H})$.
We say that $X$ is {\it maximal singularity free} (resp.\ {\it log maximal singularity free}) if $X$ does not have a maximal singularity (resp.\ log maximal singularity). 
\end{Def}

The Noether-Fano-Iskovskikh inequality \cite[Theorem 4.2]{Co95} shows that if $X$ is maximal singularity free then it is birationally superrigid.
We use (log) maximal singularity freeness to prove stability in the following sections. 
The following result shows that maximal singularity freeness characterizes birational superrigidity. 

\begin{Thm}[{Cheltsov-Shramov \cite[Theorem 1.26]{CS08}}]\label{b.rig.equivalence}
$X$ is birationally superrigid if and only if it is maximal singularity free.
\end{Thm}

We now review some examples. 
The following are Fano threefolds which have been proved to be birationally superrigid. 
\begin{itemize}
\item A smooth quartic threefold \cite{IM71}.
\item A sextic double solid, that is, a double cover $X \to \mathbb{P}^3$ ramified along a surface $S \subset \mathbb{P}^3$ of degree $6$  \cite{I79}.
\end{itemize}
The following are higher-dimensional examples.
\begin{itemize}
\item A smooth hypersurface $X_{n+1} \subset \mathbb{P}^{n+1}$ of degree $n+1$, with $n \ge 4$ \cite{dF11}, \cite{P98}.
\item A general complete intersection $X_{d_1,\dots,d_k} \subset \mathbb{P}^{n + k}$ of hypersurfaces of degree $d_i$ with $d_i \ge 2$, $\sum_{i=1}^k d_i = n + k > 3 k$ and $n \ge 4$ \cite{P01}.
\item A smooth complete intersection $X_{2,4} \subset \mathbb{P}^6$ of a quadric and a quartic which does not contain a plane \cite{C03}.
\item A double cover $X \to \mathbb{P}^n$ ramified along a hypersurface $F \subset \mathbb{P}^n$ of degree $2 n$, with $n \ge 4$ \cite{P97}.
\item A cyclic triple cover $X \to \mathbb{P}^{2 n}$ ramified along a hypersurface $F \subset \mathbb{P}^{2 n}$ of degree $3 n$ with $n \ge 2$ \cite{C04}.
\item A general cyclic cover $X \to V \subset \mathbb{P}^n$ of degree $d \ge 2$ ramified along a smooth divisor $R \subset V$ such that $V$ is a hypersurface of degree $m \ge 2$, $m + (d-1)k = n$, where $k$ is a positive integer such that $\mathcal{O}_V (R) \cong \mathcal{O}_V (d k)$, $n \ge 5$ and either $d = 2$ or $n \ge 6$ \cite{P00}, \cite{P09}. 
\item A general weighted complete intersection in a weighted projective space
\[
\mathbb{P} (1^{l+1},a_1,\dots,a_m) = \operatorname{Proj} (\mathbb{C} [x_0,\dots,x_l,y_1,\dots,y_m])
\] 
of $m+k$ hypersurfaces $y_i^2 = g_i (x_0,\dots,x_l)$, $i = 1,\dots,m$, and $f_j (x_0,\dots,x_l) = 0$, $j = 1,\dots,k$, of degree $2 l_i$ and $d_j$, respectively, such that
\[
\sum_{i=1}^m a_i + \sum_{i=1}^k d_i = l, l > 3 k \text{ and } l - k \ge 4.
\]
This is an iterated double cover of general complete intersection in 
projective space \cite{P03}. 
\end{itemize}

For each of the above examples of Fano varieties, we assume the variety is smooth. In some examples we can allow some mild singularities or drop the generality assumptions while keeping the property of birational superrigidity. 
%In all the above examples, we assume the smoothness of the respective Fano variety.
%In some examples, we can allow some mild singularities while keeping the property of %birational superrigidity and we can also drop the generality assumptions.
We refer the readers to \cite{C05} for a detailed account of this subject.
We see that every birationally superrigid Fano manifold in the above examples has index $1$ and has a base point free anticanonical divisor.  

\subsection{Seshadri constants}
Let $I\subset \mathcal{O}_X$ be a coherent ideal on $X$.
The Seshadri constant of $I$ with respect to an ample $\mathbb{Q}$-line bundle $L$ 
is defined by 
$$
	\mathrm{Sesh}(I;(X, L)):=\sup \{c\in \mathbb{Q}_{>0}\mid \pi^*L(-cE)\,\,\mbox{is ample}\},
$$
where $\pi\colon Bl_{I}(X)\to X$ is the blow up of $X$ along $I$ 
and $E$ is the associated exceptional Cartier divisor, i.e., $\mathcal{O}(-E)=\pi^{-1}I$. 
This invariant plays a key role in this paper as in 
\cite{HKLP11}, \cite{OS10}, \cite{F11} and \cite{F13}. 

\subsection{Slope stability}\label{subsec:stab}
Consult \cite[Chapter 2, especially 2.3]{Don02}, \cite[especially Section 3]{RT07} 
or \cite[Definition 2.4]{Od09} for more general background. We remark that our formulation below 
is formally different from the original presentation by Ross-Thomas \cite{RT07}, but they are equivalent, 
as proved in \cite[Theorem 4.18]{RT07}. See below for a more detailed explanation. 
Let $(X,L)$ be an $n$-dimensional polarized variety. 

A \textit{test configuration} (resp.\ a \textit{semi-test configuration}) for $(X,L)$ is a polarized scheme $(\mathcal{X},\mathcal{M})$ with a $\mathbb{G}_m$-action on $(\mathcal{X},\mathcal{M})$ and a proper flat morphism $\Pi\colon \mathcal{X}\to\mathbb{A}^1$ such that 
\begin{enumerate}
\item $\Pi$ is $\mathbb{G}_m$-equivariant for the multiplicative action of $\mathbb{G}_m$ on $\mathbb{A}^1$, 
\item $\mathcal{M}$ is relatively ample (resp.\ relatively semi-ample), and
\item $(\mathcal{X},\mathcal{M})|_{\Pi^{-1}(\mathbb{A}^1 \setminus \{0\})}$ is $\mathbb{G}_m$-equivariantly isomorphic to  $(X,L^{\otimes r})\times (\mathbb{A}^1 \setminus \{0\})$ for some positive integer $r$. 
\end{enumerate}
If $\mathcal{X}\simeq X\times \mathbb{A}^1$, $(\mathcal{X},\mathcal{M})$ is called a \textit{product test configuration}. 
Moreover, if $\mathbb{G}_m$ acts trivially, we call it a trivial test configuration. 

Slope stability deals with certain special semi-test configurations, called 
\textit{deformation to the normal cone}. The definition is as follows. Take a coherent ideal $I\subset \mathcal{O}_X$ and set $\mathcal{J}:=I+(t)\subset \mathcal{O}_{X\times \mathbb{A}^{1}}$. Then, for $r\in \mathbb{Z}_{>0}$ with 
$r > (resp.\ \geq) \frac{1}{\Sesh(I; (X,L))}$, we set $f \colon \mathcal{B}:=Bl_{\mathcal{J}}(X\times \mathbb{A}^{1})\to X\times \mathbb{A}^{1}$,  $\mathcal{L}:=f^{*}(L\times \mathbb{A}^{1})$ and $\mathcal{O}_{\mathcal{B}}(-E)=f^{-1}\mathcal{J}$ for an 
effective exceptional Cartier divisor $E$. We note that $(\mathcal{B},\mathcal{L}^{\otimes r}(-E))$ 
naturally becomes a test configuration (resp.\ semi-test configuration, if 
$\mathcal{L}^{\otimes r}(-E)$ is semiample). We call such a test configuration a \textit{deformation to the normal cone} as in \cite{RT07}. 

First, let us recall the general definition of the Donaldson-Futaki invariant 
of a test configuration $(\mathcal{X},\mathcal{M})$. 
For $k$ large, let $P(k):=\dim H^0(X,L^{\otimes k})$, which is a polynomial in $k$ of degree $n$ due to  the Riemann-Roch theorem. 
Since the $\mathbb{G}_m$-action preserves the central fibre $\mathcal{X}_0$ of $\mathcal{X}$, $\mathbb{G}_m$ acts also on $H^0(\mathcal{X}_0,\mathcal{M}^{\otimes s}\mid_{\mathcal{X}_0})$, where $s\in \mathbb{Z}_{>0}$. 
Let $w(rs)$ be the weight of the induced action on the highest exterior power of $H^0(\mathcal{X}_0,\mathcal{M}^{\otimes s}\mid_{\mathcal{X}_0})$, which is a polynomial in $s$ of degree $n+1$ due to Mumford's droll Lemma (cf.\ \cite[Lemma 2.14]{Mum77} and \cite[Lemma 3.3]{Od09}) and the Riemann-Roch theorem.
%Then, we have
%$$
%	w(k)/kP(k)=F_0+F_1k^{-1}+\mathcal{O}(k^{-2}).
%$$
Recall that, 
since an arbitrary $\mathbb{G}_m$-action on a finite-dimensional vector space is  diagonalizable, 
we can define the \textit{weights} as the exponents of the eigenvalues of the action. 
The \textit{total weight} of such an action 
%of $\mathbb{G}_{m}$ on some finite-dimensional vector space 
is defined as the sum of all  those weights. 
%We denote the total weight of the induced action on $(\Pi_{*}\mathcal{L}^{\otimes K})|_{\{0\}}$ 
%%as
%by $w(Kr)$.
% and $\dim X$ as $n$. 
%It is a polynomial of $K$ of degree $n+1$. 
%We write $P(k):=\dim H^{0}(X,L^{\otimes{k}})$. 
Now take the $rP(r)$-th power of the action of $\mathbb{G}_{m}$ on 
$(\Pi_{*}\mathcal{M})|_{\{0\}}$ and SL-normalize it. 
Then the corresponding normalized weight on 
$(\Pi_{*}\mathcal{M}^{\otimes K})|_{\{0\}}$ 
is 
$$\tilde{w}_{r,Kr}:=w(k)rP(r)-w(r)kP(k),$$ 
where $k:=Kr$. $\tilde{w}_{r,Kr}$ is a polynomial of form 
$\sum_{i=0}^{n+1}e_{i}(r)k^{i}$ of degree $n+1$ in $k$ for $k \gg 0$, with coefficients which are also 
polynomials 
of degree $\le n+1$ in $r$ 
for $r \gg 0$ : $e_{i}(r)=\sum_{j=0}^{n+1}e_{i,j}r^{j}$ for $r \gg 0$. 
Since the weight is normalized, $e_{n+1,n+1}=0$. 
The coefficient $e_{n+1,n}$ is 
called the \textit{Donaldson-Futaki invariant} of the test configuration, which we denote by $\DF(\mathcal{X},\mathcal{M})$. 
For an arbitrary \textit{semi} test configuration $(\mathcal{X},\mathcal{M})$ 
of order $r$, we can define 
%the (normalized) Chow weight or 
the Donaldson-Futaki invariant by setting $w(Kr)$ as the total weight of the  induced action on $H^{0}(\mathcal{X}, \mathcal{M}^{\otimes K})/tH^{0}(\mathcal{X}, \mathcal{M}^{\otimes K})$ (cf. \cite{RT07}). Now we can define the stability notions we are concerned with. 

\begin{Def}[{\cite{RT07}}]\label{slope.st.def}
We say that $(X,L)$ is \textit{slope stable} (resp.\ \textit{slope semistable}) 
if and only if the Donaldson-Futaki invariant is positive (resp., non-negative) 
for any non-trivial deformation to the normal cone. 
\end{Def}

Note that the above definition is a priori different from the one given in 
\cite{RT07}. Nevertheless, they are proved to be equivalent 
in \cite[proof of Theorem 4.18]{RT07}. 
(Of course, there are no 
essential differences but we will follow our formulation 
\ref{slope.st.def} just because we are more accustomed to treat the stability 
in this way.)

The following definition follows the formulation of Donaldson \cite{Don02}, 
slightly modified
\footnote{
%It is because of a subtle technical issue, that is, 
%a class of ``pathological" test configurations are found for any polarized variety: 
%whose normalizations are 
%trivial test configurations and the Donaldson-Futaki invariants vanish. 
This modification is due to a technical issue. For any polarized variety we can find a class of test configurations which are $\mathbb{G}_m$-equivariantly 
isomorphic to trivial test configuration, away from closed subschemes of codimension 
at least two. Note that those normalizations are the trivial test configurations. 
The Donaldson-Futaki invariant of the normalization of these ``pathological" test configurations automatically vanish. 
The modification of the definition excluded those class by 
considering only normal test configurations. } 
 in \cite{LX11} (cf. also \cite{Od12}). 

\begin{Def}
Suppose for simplicity that $X$ is normal. 
We say that $(X,L)$ is \textit{K-stable} (resp.\ \textit{K-semistable}) if and only if 
the Donaldson-Futaki invariant is positive (resp., non-negative) 
for any non-trivial normal test configuration. 
\end{Def}

%Let us recall that the original definition of the slope stability by Ross-Thomas 
%\cite{RT07} is of the following form.

%\begin{Def}[{Ross-Thomas \cite[Definition 4.17]{RT07}}]
%$(X,L)$ is slope stable (resp.\ slope semistable) if and only if 
%$$
%\mu_{c}(I,L)< \text{(resp.,} \leq\text{) }  \mu(X)
%$$
%for all $c\in (0,\Sesh(I;(X,L)))$ and also for 
%$c=\Sesh(I;(X,L))$ if $\Sesh(I;(X,L))\in \mathbb{Q}$ and 
%the global sections of 
%$L^{\otimes k}\otimes I^{\Sesh(I;(X,L))k}$ saturate for suffieciently divisible positive %integer $k$. 
%\end{Def}
%\noindent
%Note that the slope ``$\mu$" s, which we refer to \cite{RT07} for the precise %definitions, are defined in terms of intersection numbers on $X$ and 
%$B:=Bl_{I}(X)$. Our definition \ref{slope.st.def} above is proved in \cite[proof of Theorem %4.18]{RT07} to be equivalent to Ross-Thomas' definition. 

We end this subsection with a small remark on an extension of the above framework. 
If we take a test configuration (resp.\ semi-test configuration) $(\mathcal{X}, \mathcal{M})$, 
we can think of a new test configuration (resp.\ semi-test configuration)  $(\mathcal{X}, \mathcal{M}^{\otimes a})$ with $a \in \mathbb{Z}_{>0}$. 
From the definition of Donaldson-Futaki invariant above, we easily see that 
$\DF((\mathcal{X},\mathcal{M}^{\otimes a}))=a^{n}\DF((\mathcal{X},\mathcal{M}))$. 
%Therefore, if we consider a pair of $\mathcal{X}$ and a $\mathbb{Q}$-line bundle 
%on it, whose certain power with some linearization and some action of $\mathbb{G}_{m}$ yields a genuine test configuration (resp.\ semi-test configuration), then we can naturally define its Donaldson-Futaki invariant which is compatible with the definition above for genuine semi-test configurations. 
%We call such a pair $(\mathcal{X},\mathcal{L})$, a \textit{$\mathbb{Q}$-test configuration} (resp.\ \textit{$\mathbb{Q}$-semi-test 
%configuration}). 
%From the remark above, 
%%especially 
%in particular,
Therefore, we can define  K-stability (%resp.\
also K-polystability and K-semistability) of a pair $(X,L)$ of a projective scheme $X$ and an ample \textbf{$\mathbb{Q}$-line} bundle $L$. 
%Furthermore, our formula \ref{DF.formula} also works in this setting \textcolor{green}{so we believe this extension of the framework will often simplify arguments, though it is unnecessary and not essential.}

A key ingredient for our study is the following formula, proved in \cite{Od09}, to estimate the Donaldson-Futaki invariant for a deformation to the normal cone $(\mathcal{B}, \mathcal{L}(-E))$ derived from the (flag) ideal $\mathcal{J}:=I+(t)\subset \mathcal{O}_{X\times \mathbb{A}^1}$. 

In the following, let $X$ be a Fano $n$-fold, $I \subset \mathcal{O}_X$ an ideal and $L = \mathcal{O}_X (-r K_X)$ for some $r \in \mathbb{Z}_{>0}$ with $r \ge 1 / \Sesh (I; (X,-K_X))$.
We assume that $\mathcal{L} (-E)$ is semi-ample so that the corresponding polarized scheme $(\mathcal{B}, \mathcal{L} (-E))$ is a semi-test configuration (a deformation to the normal cone). 

\begin{Thm}[cf.\ {\cite[Theorem $3.2$]{Od09}}]\label{DF.formula}
Assume that $X$ is a Fano $n$-fold 
and let $(\overline{\mathcal{B}}:=Bl_{\mathcal{J}}(X\times \mathbb{P}^{1}),\bar{\mathcal{L}}(-E))$ be the natural compactification of $(\mathcal{B}, \mathcal{L} (-E))$.
Let $p_{i}$ $(i=1,2)$ be the projection from $X\times\mathbb{P}^{1}$ to the $i$-th factor. 
Suppose that \textcolor{black}{$\mathcal{L}(-E)$ on $\mathcal{B}$} is semi-ample.
We denote the normalization of $\mathcal{B}$ by $\tilde{\mathcal{B}}$ and use 
the same symbols for the pullbacks to $\tilde{\mathcal{B}}$ of the original polarization $\bar{\mathcal{L}}$ and  the Cartier divisor $E$. Then, the corresponding Donaldson-Futaki invariant has the following lower bound:
\begin{eqnarray}\label{10}
	\nonumber
	&&
		2n!(n+1)!\DF(\mathcal{B}, \mathcal{L}(-E))
	\\
	\nonumber
	&&\quad
		\geq
		-\big((\overline{\mathcal{L}}-E)^n.\overline{\mathcal{L}}+nE\big)
		+
		(n+1)r
		\big((\overline{\mathcal{L}}-E)^n.K_{\tilde{\mathcal{B}}/X\times \mathbb{A}^1}\big)
	\\
	\label{eq:decomposition}
	&&\quad
		=
		-\big((\overline{\mathcal{L}}-E)^n.\overline{\mathcal{L}}\big)
		+
		\big((\overline{\mathcal{L}}-E)^n.
		((n+1)rK_{\tilde{\mathcal{B}}/X\times \mathbb{A}^1}-nE)\big). 
\end{eqnarray}
\end{Thm}

\noindent
The right hand side is just the Donaldson-Futaki invariant of 
$(\tilde{\mathcal{B}},\mathcal{L}(-E))$ by \cite[Theorem3.2]{Od09} 
so that the inequality 
follows from \cite[Proposition 5.1, Remark 5.2]{RT07}. See \cite{RT07} and 
\cite{Od09} for more general statements. 

We note that $\frac{1}{r}\le \mathrm{Sesh}(I, (X,-K_{X}))$ follows from the assumption of 
the semiampleness of $\mathcal{L}(-E)$ on $\mathcal{B}$. 

\begin{Prop}[{\cite[Proposition4.3]{OS10}}]\label{prop:positivity_first_term}
The inequality 
\[
-\big((\overline{\mathcal{L}}-E)^n.\overline{\mathcal{L}}\big) \ge 0
\]
holds for any ideal $I \subset \mathcal{O}_{X}$ and equality holds if and only if $\dim (\mathrm{Supp}(\mathcal{O}_{X}/I))=0$. %change
\end{Prop}

Combining Theorem \ref{DF.formula} and Proposition \ref{prop:positivity_first_term}, we have the following criterion for the positivity of the Donaldson-Futaki invariant.

\begin{Lem}\label{q-loc}
Let $X$ be a Fano manifold and let $r \in \mathbb{Z}_{>0}$, $L = \mathcal{O}_X (-r K_X)$, $I \subset \mathcal{O}_X$ be as above.
Assume that the Weil divisor 
$$
(n+1)K_{\tilde{\mathcal{B}}/X\times \mathbb{A}^1}-n
\Sesh(I;(X,-K_{X}))E
$$ 
is effective.
Assume moreover it is nonzero if $\dim \Supp (\mathcal{O}_X/I)) = 0$.
Then $\DF (\mathcal{B},\mathcal{L} (-E)) > 0$.
\end{Lem}

Following the convention of \cite{OS10}, 
we define $a_i$ and $c_i$ as 
\begin{eqnarray*}
		K_{\tilde{\mathcal{B}}/X\times \mathbb{A}^1}
	&=&
		\sum_i a_iE_i,
	\\
		\tilde{\Pi}^{-1}\mathcal{J}
	&=&
		\mathcal{O}_{\tilde{\mathcal{B}}}(-\sum_i c_i E_i),  
\end{eqnarray*}
where $E_i$'s are exceptional prime divisors. Denote $\sum c_i E_i$ by $E$. 
The assumption of Lemma \ref{q-loc} can be reformulated as 
%follows. 
%Then, the coefficient of each exceptional prime divisor $E_i$ of the divisor 
the effectivity of 
$K_{\tilde{\mathcal{B}}/X \times \mathbb{A}^1} - n \Sesh (I;(X,-K_X)) E$ 
%is non-%negative, 
i.e., 
$$
\Sesh(I;(X,-K_{X}))
\leq \biggl(\frac{n+1}{n}\biggr)
\min_{i}
\biggl\{\frac{a_{i}}{c_{i}}\biggr\}.
$$ 
Moreover, if $\dim \mathrm{Supp}(\mathcal{O}_{X\times \mathbb{A}^1}/\mathcal{J})=0$, 
the divisor needs to be nonzero and effective, i.e., 
\begin{equation}\label{Seshineq}
\Sesh(I;(X,-K_X))<
\frac{(n+1)a_{i}}{nc_{i}} 
\end{equation}
for some $i$. 

We shall prove the main theorems in Sections \ref{div.sec} and \ref{h.cod.sec} by verifying assumptions of Lemma \ref{q-loc}.
Before going to the proof, 
we deform $I$ slightly as follows under the assumption that $|-K_X|$ is base point free. 
We write $I=\mathcal{O}_X(-F)I'$ 
where $F\in |-mK_X|$ for some $m\in \mathbb{Z}_{>0}$, 
and $I'$ is an ideal with $\codim(\Supp(\mathcal{O}/I'))\geq 2$. 
Then, we replace $F$ with a general member of $|-mK_X|$ so that 
we obtain new $I$ with the following property. 

\begin{Cond}\label{waketa}
For the blow up $\Pi \colon \mathcal{B}\rightarrow X\times \mathbb{A}^{1}$ along $\mathcal{J} := I + (t)$, if $\codim(\Pi(E_{i})\subset X\times \{0\})\geq 2$ for a $\Pi$-exceptional divisor $E_i$, then $\Pi(E_{i})\nsubseteq \Supp(F)$. 
\end{Cond}

This replacement is possible since $|-mK_X|$ is base point free by our assumption. 
Note that the Seshadri constant $\Sesh(I;(X,-K_X))$ does not change by this  
deformation. Thus, we can take the corresponding semi-test configuration 
$\mathcal{B}:=Bl_{I+(t)}(X\times \mathbb{A}^{1})$ for that perturbed $I$. 
It has the same value for the Donaldson-Futaki invariant, which follows from the description via slope (\cite{RT07}). 

Thus, we can assume Condition \ref{waketa} for $I$ from now on, 
in order to estimate the Donaldson-Futaki invariants. 

\section{Exceptional divisors with divisorial center}\label{div.sec}
In this section, under the assumptions on $X$ as in Theorem \ref{str.main}, we will prove the inequality \eqref{Seshineq} for $E_i$ in the case where $\Pi(E_i)$ is a 
divisor in $X\times \{0\}$. Let us recall that we denoted 
$\mathcal{J}=I+(t)$. As in the previous section, we write $I=I'\mathcal{O}_X (-F)$ for a coherent ideal $I' \subset \mathcal{O}_{X}$ satisfying $\codim(\Supp(\mathcal{O}_X/I'))\geq 2$ and 
a divisor $F=\sum_{j}d_{j}D_{j}$, where each $D_{j}$ is a prime divisor. 
Let $D_i$ be the component of $F$ such that $D_i \times \{0\}$ is the center $\Pi (E_i)$ of $E_i$. 

First, we claim that the inequality 
$$
\Sesh(I,(X,-K_X))\leq \frac{1}{d_i} 
$$
holds.
Let $\pi \colon X' \to X$ be the blow up of $X$ along $I'$ and write $\pi^{-1} I' = \mathcal{O}_{X'} (-E')$.
Then we have $\pi^{-1} I = \mathcal{O}_{X'} (-E'-F')$, where $F' = \pi^*F$.
Since $-K_X$ is a positive generator of $\Pic (X)$, we see that $F = \sum_j d_j D_j$ is linearly equivalent to $r (-K_X)$ for some positive integer $r$.
Note that $r \ge d_i$ since $F \geq d_i D_i$.
If the divisor
$$
\pi^*(-K_X) - c (E'+F') \sim (1 - c r) \pi^*(-K_X) - c E'
$$ is ample for some $c \in \mathbb{Q}_{>0}$ then $1 - c r > 0$.
This shows that $c < 1/r \leq 1/d_i$, hence we get the desired inequality, 
from the definition of Seshadri constant. 
%Recall that the first inequality holds in a stroger form that 
%$\Sesh(\mathcal{J};(X\times \mathbb{A}^{1},-K_{X\times \mathbb{A}^{1}}))
%=\min_{i} \{\Sesh(I_i,(X,-K_X))\}$ (\cite[Corollary 5.8]{RT07}). 

Second we have 
$$
\frac{1}{d_i}\leq \frac{a_i}{c_i}
< \frac{(n+1)a_i}{nc_i}, 
$$ 
since the pair $(X\times \mathbb{A}^{1}, D_{i}\times \mathbb{A}^{1})$ 
is canonical around the generic point of $D_i\times \{0\}$ and the discrepancy of 
$(X \times \mathbb{A}^1, D_i \times \mathbb{A}^1)$ at $E_i$ is $a_i - c_i/d_i$. 
As a conclusion, we get the desired inequality \eqref{Seshineq} for $E_i$ with $\dim(\Pi(E_i))=n-1$.
We note that the condition $\rho=\operatorname{index}(X)=1$ is sufficient for the arguments in this section. 

\section{Exceptional divisors with higher codimensional center}\label{h.cod.sec}
In this section, under the same assumptions on $X$, 
we will prove the inequality \eqref{Seshineq} for $E_i$ in the case where 
$\codim(\Pi(E_i)\subset X\times \{0\})\geq 2$. 
Recall that we have Condition \ref{waketa} for $I$. 

We first note that $\Sesh (I, (X,-K_X)) \le \Sesh (I',(X,-K_X))$ since $\Supp (\mathcal{O}_X/I')$ is not contained in $\Supp (F)$.
Take a positive rational number $c< \Sesh(I,(X,-K_{X})) \le \Sesh (I',(X,-K_X))$. 
For a sufficiently divisible positive interger $l$, 
set the linear system $\Sigma^{(c)}_{I',l}\subset |-lK_{X}|$ which corresponds to 
$H^{0}(X,{I'}^{cl}\mathcal{O}_{X}(-lK_X))\subset H^{0}(X,\mathcal{O}_{X}(-lK_X))$. 
The linear system $\Sigma^{(c)}_{I',l}$ is movable since $\codim \Supp (\mathcal{O}_X/I') \ge 2$ and $c < \Sesh (I', (X,-K_X))$.
By log maximal singularity freeness, 
the pair $(X, \frac{1}{l} \Sigma^{(c)}_{I',l})$ is log canonical, hence the pair $(X,\frac{1}{l}D)$ is log canonical for a general member $D$ of $\Sigma^{(c)}_{I',l}$.
By the inversion of adjunction on log canonicity \cite[17.7 Theorem]{K+92}, the pair $(X \times \mathbb{A}^1, \frac{1}{l} (D \times \mathbb{A}^1) + X \times \{0\})$ is log canonial near $X \times \{0\}$.
The discrepancy of the latter pair at $E_i$ is $a_i - \frac{1}{l} l c c'_i - b_i = a_i - c c'_i - b_i$, where $c_{i}':=\val_{E_{i}}(I')$, the algebraic valuation of $I'$ with respect to $E_i$, and $b_i$ is defined by $\tilde{\Pi}^* (X \times \{0\}) = \tilde{\Pi}_*^{-1} (X \times \{0\}) + \sum_i b_i E_i$.
Note that $b_i \ge 1$ since $X \times \{0\}$ is a Cartier divisor and the center of $E_i$ on $X \times \mathbb{A}^1$ is contained in $X \times \{0\}$.
The log canonicity implies $a_i - c c'_i - b_i \ge -1$, which together with $b_i \ge 1$ shows the inequality $c \le a_i/c'_i$.
Since this holds for any such $c$, we have 
\begin{equation}\label{I_0'}
\Sesh(I,(X,-K_{X}))\leq \min_{E_{i}\subset \Exc(\Pi)}\biggl\{\frac{a_{i}}{c_{i}'}\biggr\}, 
\end{equation}
where $\Exc(\Pi)$ is the exceptional locus of $\Pi$. 
If $F\equiv -mK_X$ with $m\in \mathbb{Z}_{>0}$, 
the inequality (\ref{I_0'}) is equivalent to the inequality 
\begin{equation}\label{I_0}
\Sesh(I,(X,-K_{X}))\leq 
\dfrac{\min\{\frac{a_i}{c_i'}\}}
{1+m\cdot\min \{\frac{a_i}{c_i'}\}}. 
\end{equation}
By Condition \ref{waketa}, 
we have $c_{i}=\val_{E_{i}}(\mathcal{J})\leq \val_{E_{i}}(I)=\val_{E_i}(I')=:c_i'$. 

Summing up, 
\begin{equation}\label{highcodim}
\Sesh(I,(X,-K_{X}))
\leq 
\dfrac{\min\{\frac{a_i}{c_i'}\}}
{1+m\cdot\min \{\frac{a_i}{c_i'}\}}
\leq
\frac{a_i}{c_i'}
\leq
\frac{a_i}{c_i}
<\frac{(n+1)a_i}{nc_i}. 
\end{equation}
We note that the replacement of the divisorial part $F$ of $I = \mathcal{O}_X (-F) I'$ 
(which makes use of the assumption of the base point freeness of $|-K_X|$) 
gives the last inequality. Therefore, we have the desired inequality and
the proof of Theorem \ref{str.main} is completed by Lemma \ref{q-loc}. 
$\square$\\

Now we have completed the proof of Theorem \ref{str.main}. 
Theorem \ref{thm:main} follows immediately from Theorem \ref{str.main} since log maximal singularity freeness implies maximal singularity freeness. 

\begin{Rem}
A difficulty arises when we try to strengthen Theorems \ref{thm:main} and 
\ref{str.main} to K-stability. 
The point is that a flag ideal $\mathcal{J}=I_0+I_1+\cdots I_{N-1}t^{N-1}+(t^N)$  of length $N>1$ should satisfy a condition that $I_i\subset I_j$ for any $i<j$, 
which makes it hard to deform it to an another flag ideal which satisfies Condition \ref{waketa}. 
\end{Rem}

\begin{Rem}
We explain a generalization of the main theorems.
Let $X$ be a Fano manifold and $G \subset \Aut (X)$ a finite subgroup.
One can define $G$-{\it birational superrigidity} (see \cite[Definition 1.30]{CS08} for the definition).
On the other hand, one can also define $G$-{\it equivariant slope stability} 
as the condition that $\DF (\mathcal{B}, \mathcal{L} (-E)) > 0$ for every deformation to the normal cone derived from the ideal $\mathcal{J} = I + (t) \subset \mathcal{O}_{X \times \mathbb{A}^1}$, where $I$ is a $G$-invariant ideal of $\mathcal{O}_X$ (cf. \cite[Section 2.2]{OS10}). 
Note that when $G = \{ \operatorname{id} \}$ is the trivial group the $G$-birational superrigidity (resp., $G$-equivariant slope stability) coincides with usual birational superrigidity (resp., slope stability). 

We can prove the following result with minor natural modifications.
Let $X$ be a Fano manifold.
Assume that there is a finite subgroup $G \subset \Aut (X)$ with the following properties.
\begin{enumerate}
\item $X$ is $G$-birationally superrigid.
\item For every $G$-invariant nonzero effective divisor $D$ on $X$, one has $D \ge -K_X$.
\item Every $G$-invariant nonzero effective divisor $D$ on $X$ is $G$-base point free.
\end{enumerate}
Then $X$ is $G$-equivariantly slope stable. 

Here, an effective divisor $D$ on $X$ is $G$-{\it base point free} if for any point $x \in X$ there is a $G$-invariant effective divisor $D' \in |D|$ which does not contain $x$.

We explain how to modify our arguments in order to prove the above generalization.
It is enough to prove the inequality \eqref{Seshineq} assuming that $I \subset \mathcal{O}_X$ is $G$-invariant.
By the condition (iii) above, we can deform $I$ keeping the $G$-invariance so that Condition \ref{waketa} holds, which we assume without loss of generality from now on. 
In order to prove inequality \eqref{Seshineq} for exceptional divisors with divisorial centers, it is enough to prove the inequality $\Sesh (I, (X, -K_X)) \le 1/d_i$ (see Section \ref{div.sec} for $d_i$), which follows from condition (ii) above.
In order to prove inequality \eqref{Seshineq} for exceptional divisors with higher codimensional centers, it is enough to show that $\lct (X, \Sigma) \geq 1/l$, where $\Sigma := \Sigma^{(c)}_{I',l} \subset |- l K_X|$ is a $G$-invariant movable linear system constructed as in Section \ref{h.cod.sec}.
It follows from condition (i) above and \cite[Theorem 1.31]{CS08} that $\lct (X, \Sigma) \geq \ct (X, \Sigma) \ge 1/l$, which completes the proof.
\end{Rem}

\section{A conjecture}

Note that we treat in this paper a special class of Fano manifolds of Picard rank $1$ (i.e., $\Pic(X)\cong\mathbb{Z}$). 
More generally, we expect the following. 

\begin{Conj}
For an arbitrary Fano manifold $X$ of Picard rank $1$, 
$(X,-K_X)$ is K-semistable. 
\end{Conj}

We note some supporting evidences here. 
First, it is proved in \cite{F13} that $(X,-K_X)$ is slope stable with respect to divisors. 
Recall also that we proved a stronger statement in section \ref{div.sec} under the additional assumption of Fano index $1$. Moreover, the results of Hwang, Kim, Lee, Park \cite[Theorem1.3]{HKLP11} and Fujita \cite[Theorem 1.1]{F11}, also show slope semistablity along smooth curves. 

We also remark that we cannot expect strict (poly)stability as, for instance, 
small deformations of the Mukai-Umemura $3$-fold are \textit{not} K-polystable but K-semistable (see \cite{Tia97}).

\end{document}